\documentstyle[11pt]{article}

\input{psfig}

\begin{document}
\vspace*{.6in}

\centerline {{\Large {\bf Dynamics of Certain Smooth One-dimensional
Mappings}}}

\vskip10pt
\centerline {\Large \bf II. Geometrically finite
one-dimensional mappings}

\vskip25pt
\large
\centerline{Yunping Jiang }
\centerline{Institute for Mathematical Sciences, SUNY at Stony Brook}
\centerline{Stony Brook, L.I., NY 11794}
\vskip20pt
\centerline{September 1, 1990}

\vskip40pt

\centerline{ {\Large \bf Abstract}}

\vskip10pt

We study geometrically finite one-dimensional mappings.
These are a subspace of $C^{1+\alpha
}$ one-dimensional mappings with finitely many, critically
finite critical
points.  We study some geometric properties of a mapping in this
subspace. We prove that this subspace is closed under
quasisymmetrical conjugacy.  We also prove that if two mappings in
this subspace are topologically conjugate, they are
then quasisymmetrically conjugate. We show some examples of geometrically
finite one-dimensional mappings.

\pagebreak
\centerline{{\Large \bf Contents}}

\vskip20pt
\noindent {\Large \S 1 Introduction.}

\vskip10pt
\noindent {\Large  \S 2  Geometrically Finite One-dimensional Mappings.}

\vskip5pt
\noindent \S 2.1 Geometrically finite.

\vskip5pt
\noindent \S 2.2 Bounded geometry.

\vskip5pt
\noindent \S 2.3 From geometrically finite to bounded geometry.

\vskip5pt
\noindent \S 2.4 Quasisymmetrical classification.

\vskip5pt
\noindent \S 2.4.1 Bounded nearby geometry.

\vskip5pt
\noindent \S 2.4.2 Quasisymmetry.

\vskip5pt
\noindent \S 2.5 Closeness under quasisymmetrical conjugacy.

\vskip10pt
\noindent {\Large \S 3
Examples Of Geometrically Finite One-dimensional Mappings.}

\vskip5pt
\noindent \S 3.1 A $C^{3}$-mapping with nonpositive Schwarzian derivative.

\vskip5pt
\noindent \S 3.2 A $C^{1}$-mapping with bounded variation derivative.

\vskip5pt
\noindent \S 3.3 A question on $C^{1+\alpha}$-mappings with expanding
periodic points.

\pagebreak

\centerline {\Large {\bf \S 1 Introduction}}

\vskip10pt

{\bf Quasisymmetrical conjugacy.} Two smooth mappings $f$ and $g$ from a
one-dimensional manifold $M$ to itself are
topologically conjugate if there is a homeomorphism $h$ from $M$ to
itself such that $f\circ h =h\circ g$. The homeomorphism $h$ and its inverse
are usually
not both Lipschitz; if they are, then all the eigenvalues of $f$ and $g$ at the
periodic points have to be the same. Between the class of homeomorphisms and the class of
Lipschitz homeomorphisms, there is a class of \underline {quasisymmetric homeomorphisms}.
A quasisymmetric
homeomorphism distorts symmetrically placed triples by a
bounded amount.
A celebrated Ahlfors-Beurling extension theorem [A] tells us that any quasisymmetric
homeomorphism of the real line can be extended to
a quasiconformal homeomorphism of the complex plane.
Thus quasisymmetric
property of the conjugating homeomorphism gives us a chance to use some methods and
theorems in one complex variable functions to study the dynamics of some smooth
one-dimensional
mappings.
M. Jakobson recently
considered
a $C^{3}$-folding mapping with negative Schwarzian derivative
and one non-recurrent critical point. He proved that if two such
mappings are topologically conjugate, they are then quasisymmetrically
conjugate [Ja].  D. Sullivan [S1], M. Herman [H], J. Yoccoz [Y] and G.
Swiatek [SW], etc.,  have some interesting
results on this direction for some folding mappings and critical circle mappings.
\vskip7pt
{\bf What we would like to say in this paper.}
We consider a subspace of piecewise $C^{1+\alpha }$-mappings
with finitely many, critically finite critical points from a compact smooth
one-dimensional manifold into itself and study some geometric properties
of a mapping in this subspace.

\vskip7pt
Suppose $M$ is an oriented connected compact one-dimensional
$C^{2}$-Riemannian manifold with Riemannian metric $dx^{2}$ and associated
length element $dx$. Suppose $f: M\mapsto M$ is a $C^{1}$-mapping.
Furthermore, without loss generality, we will assume that $f$ maps the
boundary of $M$ (if it is not empty) into itself and the one-sided
derivatives of $f$
at all boundary points of $M$ are not zero.

\vskip5pt
We say $c\in M$ is a
critical point of $f$ if
the derivative of $f$ at this point is
zero.
We say a critical point of $f$ is \underline{critically
finite} if its orbit consists
of finitely many points.

\vskip5pt
Suppose $f: M\mapsto
M$ is a $C^{1}$-mapping with finitely many,
critically finite critical points.
There is a natural \underline{Markov partition} of $M$
by $f$. This Markov partition consists of the intervals of the complement
of the critical
orbits of $f$. We call it the first partition $\eta_{1}$ of $M$ by $f$.
For any positive
integer $n$, the $n^{th}$-partition $\eta_{n}$ of $M$ by $f$
consists of all the intervals $I'$ such that the restriction of the $(n-1)^{th}$-iterate of $f$ is a homeomorphism
from it to an interval in the first partition of
$M$ by $f$.  We use $\lambda_{n}$ to denote the maximum of the
lengths of the intervals in the $n^{th}$-partition of $M$ by $f$.
We say the $n^{th}$-partition of $M$ by $f$ goes to zero
\underline{exponentially with $n$} if there are constants
$K>0$ and $0< \mu<1$
such that $\lambda_{n}\leq K\mu^{n}$ for every $n$.

\vskip5pt
\underline{A geometrically finite one-dimensional mapping} is a
$C^{1+\alpha }$-mapping $f: M\mapsto M$ for some $0 < \alpha \leq 1$
with finitely many, critically finite, non-periodic power law critical
points such that
the $n^{th}$-partition
of $M$ by $f$ goes to zero
exponentially with $n$. The reader may see \S 2 for a definition of a power
law critical point of $f$. We also note that the definition of $C^{1+\alpha}$ for a
mapping with power law critical points is given in \S 2 and is little
different from the usual one.

\vskip5pt
To study a geometrically finite one-dimensional mapping, we introduce
two concepts,
\underline{bounded geometry} and \underline{bounded nearby geometry},
for
a sequence $\eta =\{ \eta_{n} \}_{n=1}^{\infty }$ of nested partitions.
We
say a sequence $\eta=\{ \eta_{n} \}_{n=1}^{\infty}$ of nested partitions
has bounded geometry
if there is a positive constant $K$ such that for any  $J\subset I$ with $J\in \eta_{n+1}$
and $I\in \eta_{n}$, the ratio of lengths,
$|J|/|I|$, is bounded by $K$ from below. We say this sequence has
bounded nearby geometry if there is a positive constant $K$ such that
for any $J$ and $I$ in $\eta_{n}$ with a common endpoint, the ratio of
lengths, $|J|/|I|$,  is bounded by $K$ from below.  The
bounded geometry here is an analogue to the Sullivan's
definition of bounded geometry for a Cantor set on the line [S2].
One of the main theorems in this paper is the following (see Theorem A and Lemma 2).

\vskip5pt
{\sc Main Theorem.} {\em  Suppose $f:M\mapsto M$ is geometrically
finite and $\eta=\{ \eta_{n}\}_{n=1}^{\infty}$ is the induced sequence
of nested partitions of $M$ by $f$.
Then
the sequence $\{ \eta_{n} \}_{n=1}^{\infty }$ of
nested partitions has bounded geometry and bounded nearby geometry.}

\vskip5pt
The proof of this theorem is an
application of the $C^{1+\alpha}$-Denjoy-Koebe distortion lemma in [J2].

\vskip5pt
Following the methods in [MT], we can classify topologically the
geometrically finite one-dimensional mappings by their kneading invariants.
Moreover, using these properties, bounded geometry and bounded nearby
geometry, we can classify these mappings quasisymmetrically as follows.

\vskip5pt
A homeomorphism $h: M\mapsto M$ is
quasisymmetrical
if there is a positive constant $K$ such that for any
two points $x$ and $y$ in $M$ and $z=(x+y)/2$,
\[ K^{-1} \leq \frac{|h(x)-h(z)|}{|h(z)-h(y)|} \leq K.\]
We say two mappings $f$ and $g$ from $M$ to itself
are quasisymmetrically
conjugate if they are topologically conjugate and the conjugating
homeomorphism is quasisymmetrical.

\vskip5pt
{\sc Theorem B.} {\em Suppose $f$ and $g$ are geometrically finite and
topologically
conjugate. They are then quasisymmetrically conjugate.}

\vskip5pt
Geometrically finite one-dimensional mappings are
closed under quasisymmetrical conjugacy in the space of
$C^{1+\alpha}$-mappings with only power law critical points as follows.

\vskip5pt
{\sc Theorem C.} {\em If a $C^{1+\alpha }$-mapping
$f:M\mapsto M$ for some $0< \alpha \leq 1$ with only power law critical
points is quasisymmetrically conjugate to a geometrically
finite one-dimensional mapping, then it is also a geometrically finite
one-dimensional mapping.}

\vskip5pt
One example of a geometrically finite
one-dimensional mapping is the following (see Section 2 for details).

\vskip5pt
{\sc Example 1.} {\em A $C^{3}$-mapping $f:M\mapsto M$ with
nonpositive Schwarzian derivative and finitely many, critically
finite, nonperiodic power law critical points.}

\vskip5pt
This kind of mappings was systematically studied by
M. Misiurewicz in 1979 [Mi] and many other people [Ja], [BL] and [MS],
etc.

\vskip5pt
Let $C^{1+bv}$ stand for $C^{1}$ with bounded
variation derivative (the definition of $C^{1+bv}$ for a
mapping with power law critical points is given in \S 3.2). We say a
periodic point $p$ of a mapping $f$ is
expanding if the absolute value of the eigenvalue (some people
call an eigenvalue an multiplier) $(f^{\circ n})^{\prime }(p)$ is
greater than one,  where $n$ is the period of $p$.
The main theorem in \S 3 is the following.

\vskip5pt
{\sc Theorem D.} {\em Suppose $f:M\mapsto M$ is a $C^{1+\alpha}$-, for some $0<
\alpha \leq 1$, and
$C^{1+bv }$-mapping with finitely many, critically finite, nonperiodic power
law critical points and only expanding periodic points
and suppose $\eta =\{ \eta_{n} \}_{n=1}^{\infty}$ is the induced
sequence of nested partitions of $M$ by $f$. Then $\eta$ has bounded
geometry.}

\vskip5pt
The study of this theorem is inspired by the paper [M] where R. Ma\~ne [M]
proved that a $C^{2}$-endomorphism $f:M\mapsto M$ with
only expanding periodic points is actually expanding
in a suitable smooth coordinate on $M$.
Theorem D provides
another example of a geometrically finite one-dimensional mapping.

\vskip5pt
{\sc Example 2.} {\em A $C^{1+\alpha}$- , for some $0< \alpha \leq 1$, and
$C^{1+bv}$-mapping $f:M\mapsto M$
with finitely many, critically finite, nonperiodic power law critical
points and only expanding periodic points.}

\vskip5pt
In Theorem D and Example 2, the condition that $f$ is a
$C^{1+\alpha}$- , for some $0< \alpha \leq 1$, and $C^{1+bv }$-mapping can not be
weakened
to the condition that $f$ is a $C^{1+\alpha }$-mapping for there is a
counterexample in [J1].
The construction
of the counterexample in [J1] is like the
construction of the Denjoy counterexample in circle diffeomorphisms and
this example is not topologically conjugate to any geometrically
finite one-dimensional mapping.

\vskip5pt
The condition that a $C^{1+\alpha }$-mapping for some $0< \alpha \leq
1$ with only power law critical points is quasisymmetrically conjugate
to a geometrically finite one-dimensional mapping in Theorem C
can not be weakened
to the condition that a $C^{1+\alpha}$-mapping with only power law critical
points is topologically
conjugate to
a geometrically finite one-dimensional mapping too for there is an easy
counterexample (see Figure 4 in \S 3.3). This counterexample
has a neutral fixed point (namely the absolute value of the
eigenvalue of $f$ at this fixed point is one) and
suggests a question
as follows.

\vskip5pt
{\sc Question 1.} {\em Suppose $f: M\mapsto M$
is a $C^{1+\alpha}$-mapping for some
$0< \alpha \leq 1$
with
only power law critical points
and only expanding periodic points
and is topologically conjugate to a geometrically
finite one-dimensional mapping.
Is $f$ geometrically finite ?}

\vskip20pt
{\bf Acknowledgment.} The author would like to thank
Dennis Sullivan and John Milnor for their constant encouragement and
many conversations. I would
also like to thank Elise Cawley, Benjamin Bielefeld, Mikhail Lyubich,
Scott Sutherland, Grzegorz Swiatek,
Folkert
Tangerman and Peter Veerman
for many useful conversations and help.

\vskip20pt
\centerline {\Large \bf  \S 2 Geometrically Finite One-dimensional Mappings}

\vskip3pt
Suppose $M$ is an oriented connected compact one-dimensional
$C^{2}$-Riemannian manifold with Riemannian metric $dx^{2}$ and associated
length element $dx$.  Suppose $f: M\mapsto M$ is a
continuous mapping. We say an interior point
$c\in M$ is a critical point if

\vskip3pt
$(a)$ $f$ is not differentiable at $c$, or

\vskip3pt
$(b)$ $f$ is differentiable
at $c$ and the derivative of $f$ at $c$ is zero.

\vskip5pt
We always assume that $f$ is
$C^{1}$ at any non-critical point $p$, namely $f$ is differentiable in a
small neighborhood $U_{p}$ of $p$ and the derivative $f'$ of
$f$ in the neighborhood $U_{p}$ is continuous. We say a
critical point $c$ is a power law critical
point if

\vskip3pt
$(c)$ $c$ is an isolated critical point,

\vskip3pt
$(d)$ for some $\gamma \geq 1$,
\[ \lim_{x\mapsto c+} \frac{f'(x)}{|x-c|^{\gamma -1}} \hskip5pt and \hskip5pt
\lim_{x\mapsto c-} \frac{f'(x)}{|x-c|^{\gamma-1}} \]
have nonzero limits $A$ and $B$.

\vskip5pt
We call the numbers $\gamma$ and $\tau =A/B$ the exponent and the
asymmetry of $f$ at $c$ (see [J2]).
We say a critical point $c$ of $f$ is critically finite if
the orbit $\{ c,$ $f(c)$, $\cdots \}$ is a finite set.

\vskip5pt
Although the results in this paper hold for a piecewise $C^{1}$-mapping
$f:M\mapsto M$ with both smooth and non-smooth critical points, but we are
only interested in a smooth critical point of $f$. Henceforth we will assume
that $f: M\mapsto M$ is a $C^{1}$-mapping.
Furthermore, without loss generality, we will assume that $f$ maps the
boundary of $M$ (if it is not empty) into itself and the one-sided
derivatives of $f$ at all boundary points of $M$
are not zero. We note that in the general case, a boundary point
of $M$ should count as a critical point anyhow.

\vskip5pt
We define the term $C^{1+\alpha }$ for a real number $0< \alpha \leq 1$.
Suppose $f:M\mapsto M$ has only power law critical points.
We use $CP=\{ c_{1}$,
$\cdots $, $c_{d} \}$ to denote the set of critical points
of $f$ and use $\Gamma=\{ \gamma_{1}$, $\cdots $, $\gamma_{d} \}$
to denote the corresponding exponents of $f$. Suppose $\eta_{0}$ is the set of
the closures of the intervals of the complement of the set of critical
points $CP$ of
$f$ in $M$.

\vskip5pt
{\sc Definition 1.} {\em We say the mapping $f$ is
$C^{1+\alpha}$ for some
$0< \alpha \leq 1$ if

\vskip3pt
$(*)$ the restrictions of $f$ to the intervals in $\eta_{0}$
are $C^{1}$ with $\alpha $-H\"older continuous derivatives
and

\vskip3pt
$(**)$ for every critical point $c_{i}$ of $f$, there
is a small neighborhood $U_{i}$ of $c_{i}$ in $M$ such that
$r_{-,i}(x)=f'(x)/|x-c|^{\gamma_{i}-1}$ for $x<c$ in $U_{i}$ and
$r_{+,i}(x)=f'(x)/|x-c|^{\gamma_{i}-1}$ for $x>c$ in $U_{i}$ are
$\alpha$-H\"older continuous functions.}

\vskip5pt
We define the term exponential decay. Suppose $f:M\mapsto M$ is a
$C^{1}$-mapping such that the set of
critical orbits $\cup_{n=0}^{\infty}f^{\circ n}(CP)$ is finite.
Let $\eta_{1}$ be the
set $\{ I_{1}$, $\cdots $, $I_{n}\}$ of the closures of the intervals of the
complement
of the critical orbits $\cup_{n=0}^{\infty}f^{\circ n}(CP)$ in $M$.
We call it the first partition of $M$ by $f$. It is a
Markov partition, namely
$f$ maps every interval in it
into and onto the union of some intervals in it.  Let $\eta_{n}=
f^{-(n-1)}(\eta_{1})$ be the set of all the intervals, to each of which
the restriction of
the $(n-1)^{th}$-iterate of $f$ is a homeomorphism from this interval
to an interval in the first
partition $\eta_{1}$. We call it the $n^{th}$-partition of $M$ by $f$.	We use $\eta $
to denote the sequence $\{ \eta_{n} \}_{n=1}^{\infty }$ of nested partitions and call it the
induced sequence of nested partitions of $M$ by $f$.
Let $\lambda_{n}$ be the
maximum of lengths of
the intervals in the $n^{th}$-partition $\eta_{n}$.
We say the $n^{th}$-partition $\eta_{n}$ tends to zero
exponentially with $n$ if
there are constants $K>0$ and $0< \mu <1$
such that $\lambda_{n}\leq K\mu^{n}$ for
all the positive integers $n$.

\vskip5pt
\noindent {\bf \S 2.1 Geometrically finite.}

\vskip3pt
We now give the definition of a geometrically finite one-dimensional
mapping as follows.

\vskip5pt
{\sc Definition 2.} {\em  We say a mapping $f: M\mapsto M$ with only
power law
critical points is geometrically finite if it satisfies the following
conditions:

\vskip3pt
Smooth condition: $f$ is $C^{1+\alpha }$ for
some $0 < \alpha \leq 1$.

\vskip3pt
Finite condition:
the set of critical orbits $\cup_{i=0}^{\infty }f^{\circ }(CP)$
is finite.

\vskip3pt
No cycle condition: no critical point is a periodic point of $f$.

\vskip3pt
Exponential decay condition: the $n^{th}$-partition
$\eta_{n}$ tends to zero exponentially with $n$.}

\vskip5pt
\noindent {\bf \S 2.2 Bounded Geometry.}

\vskip3pt
We say a set of finitely many closed subintervals of $M$ with pairwise
disjoint
interiors is a partition of $M$ if the union of these intervals is $M$.
Suppose $\eta =\{ \eta_{n}\}_{n=1}^{\infty}$ is a sequence of partitions
of $M$. We say it is nested if every interval in $\eta_{n}$ is the union
of some intervals in $\eta_{n+1}$ for every $n \geq 1$.

\vskip5pt
{\sc Definition 3.} {\em We say a sequence $\eta =\{
\eta_{n}\}_{n=1}^{\infty}$
of nested partitions has bounded geometry if there is a positive
constant $K$ such that for any pair $J\subset I$ with $J\in \eta_{n+1}$
and $I\in \eta_{n}$, the ratio
$|J|/|I| \geq K$. We call the biggest such constant $BC_{f}$ the bounded geometry constant.}

\vskip5pt
\noindent {\bf \S 2.3 From geometrically finite to bounded geometry.}

\vskip3pt
One of the main theorems in this paper is the following:

\vskip3pt
{\sc Theorem A.} {\em Suppose $f:M\mapsto M$ is geometrically finite
and $\eta=\{ \eta_{n}\}_{n=0}^{\infty}$ is the induced sequence of
nested partitions of $M$ by $f$. Then $\eta$ has bounded geometry.}

\vskip3pt
Before to prove this theorem, let me state the $C^{1+\alpha}$-Denjoy-Koebe
distortion lemma in [J2]. For a
geometrically finite one-dimensional mapping, this
lemma can be written in the following simple form (see \S 3.3 in
[J2]).

\vskip3pt
{\sc Lemma 1.} (The $C^{1+\alpha}$-Denjoy-Koebe distortion lemma)
{\em Suppose $f:M\mapsto M$ is geometrically finite.
There are two positive
constants $A$ and $B$ and a positive integer $n_{0}$ such that for
any inverse branch $g_{n}$ of $f^{\circ n}$ and
any pair $x$ and $y$ in the intersection of one of the intervals in
$\eta_{n_{0}}$ and the domain of $g_{n}$,
the distortion $|g_{n}(x)/g_{n}(y)|$ of $g_{n}$ at these two points satisfies
\[ \frac{|g_{n}(x)|}{|g_{n}(y)|} \leq \exp \Big( A+\frac{B}{D_{xy}} \Big) \]
where $D_{xy}$ is the distance between $\{ x$, $y\}$ and the
post-critical orbits $\cup_{i=1}^{\infty}f^{\circ i}(CP)$.}

\vskip5pt
{\it Proof of Theorem A.}
Suppose $n_{0}$, $A$ and
$B$ are the constants in Lemma 1. Suppose $\{ c_{i_{1}}$, $\cdots $,
$c_{i_{k}}\}$ is a sequence of critical points of $f$.
We say it is a critical chain of $f$ if there is a sequence
$\{ l_{i_{1}},$ $\cdots$, $l_{i_{k-1}}\}$ of the integers such that
$f^{\circ l_{1}}(c_{i_{1}})=c_{i_{2}}$, $\cdots $,
$f^{\circ l_{k-1}}(c_{i_{k-1}})=c_{i_{k}}$. We call the integer
$l=l_{i_{1}}+\cdots l_{i_{k-1}}$ the length of this chain.  By
the no
cycle condition, there are only finitely many critical chains. Let
$N_{0}$ be the maximum of lengths of all the critical chains of $f$.

\vskip3pt
We say an interval in $\eta_{n}$ is a critical interval if one of
its endpoints is a critical point.
We may assume that for every
critical interval in $\eta_{n_{0}}$, one of its endpoints is not in the
critical orbits $\cup_{i=0}^{\infty}
f^{\circ i}(CP)$. Let ${\cal U}$ be the union of all the critical
intervals in
$\eta_{n_{0}}$ and $K_{1}>0$ be the minimum of ratios, $|J|/|I|$, for
$J\subset I$ with $J\in \eta_{n_{0}+1}$ and $I\in \eta_{n_{0}}$.

\vskip3pt
For any $J\subset I$ with $J\in \eta_{n+1}$, $I\in \eta_{n}$ and $n> n_{0}$, let
$J_{i}=f^{\circ i}(J)$ and $I_{i}=f^{\circ i}(I)$ for $i=0$, $\cdots $,
$n-n_{0}$. Then
$J_{n-n_{0}}\in \eta_{n_{0}+1}$ and $I_{n-n_{0}}\in \eta_{n_{0}}$.  We consider the intervals $\{ I_{0},$ $\cdots $, $I_{n-n_{0}}\} $ in the two cases. One is
 that no one of them is in ${\cal U}$. The other
is that at least one of them is in ${\cal U}$.

\vskip3pt
For the first case, by using the naive distortion lemma (see [J1] or
[J2]),
there is a constant $K_{2}>0$ (which does not depend on any particular intervals $J\subset I$)
such that for any $x$ and $y$ in $I$,
\[ \frac{|f^{\circ (n-n_{0})}(x)|}{|f^{\circ (n-n_{0})}(y)|} \geq
K_{2},\]
and moreover,
\[ \frac{|J|}{|I|} \geq K_{3}=K_{2}K_{1}.\]

\vskip3pt
For the second case, let $l\leq n-n_{0}$ be the greatest integer such that
$I_{l}\subset {\cal U}$. We note that $I_{i}\cap {\cal U}=\emptyset$ for $i=l+1$, $\cdots $,
$n-n_{0}$.  By using the naive distortion lemma like that in the first case,
 we can also show that
\[ \frac{|J_{l+1}|}{|I_{l+1}|} \geq K_{3}=K_{2}K_{1}.\]
Let $\tilde{I}_{i}$ be the interval in $\eta_{n_{0}+l-i}$ containing
$I_{i}$ for $i=0$, $\cdots $, $l$.
Then $\tilde{I}_{l}$ is an interval in $\eta_{n_{0}}$ and is contained in ${\cal U}$. Suppose
$c_{q}\in CP$ is an endpoint of $\tilde{I}_{l}$. The restriction
of $f$ to $\tilde{I}_{l}$ is comparable to the mapping
$x:\mapsto |x-c_{q}|^{\gamma_{q}}+f(c_{q})$. We can find a
positive constant $K_{4}$ (only depends on $K_{3}$) such that
\[ \frac{|J_{l}|}{|I_{l}|} \geq K_{4}.\]

\vskip3pt
We may assume that both endpoints of $\tilde{I}_{l}$ are not in the
post-critical orbits $\cup_{n=1}^{\infty}f^{\circ n}(CP)$.
Otherwise, by the no
cycle condition, there is $k\leq N_{0}$ such that $\tilde{I}_{l-k}$ has
this property, one of its endpoint is a critical point of $f$
and both of its endpoints are not in the post-critical orbits
$\cup_{n=1}^{\infty}f^{\circ n}(CP)$. Then we can use
$\tilde{I}_{l-k}$ to instead of $\tilde{I}_{l}$ because there is a
constant $K_{5}$ (only depends on $K_{4}$) such that
\[ \frac{|J_{l-k}|}{|I_{l-k}|} \geq K_{5}\frac{|J_{l}|}{|I_{l}|}.\]

\vskip3pt
Let $K_{6}$ be the minimum of lengths of the critical intervals
in $\eta_{n_{0}}$. Now using Lemma 1 (the
$C^{1+\alpha }$-Denjoy-Koebe distortion lemma), for any $x$ and $y$ in
$I_{l}\subset \tilde{I}_{l}$,
\[ \frac{ |f^{\circ (n-l)}(x)|}{ |f^{\circ (n-l)}(y)|}\geq K_{7}=-\exp
\Big( A+\frac{B}{K_{6}}\Big) ,\] and moreover,
\[ \frac{|J|}{|I|} \geq C_{8}=K_{7}K_{4}.\]
The bounded geometry constant $BC_{f}$ is greater that the maximum of
$K_{3}$ and $K_{8}$.

\vskip5pt
\noindent {\bf \S 2.4 Quasisymmetrical classification.}

\vskip3pt
Topologically, we can classify the geometrically
finite one-dimensional mappings by their kneading invariants
just following the methods in [MT] (see [MT] for a definition of a
kneading sequence). By this we means
that for two geometrically finite one-dimensional mappings $f$ and $g$,
there is an orientation-preserving homeomorphism $h: M
\mapsto M$ such that $f\circ h=h\circ g$ if and only if the kneading
invariants of $f$ and $g$ are the same.

\vskip3pt
A homeomorphism $h: M\mapsto M$ is
quasisymmetrical
if there is a positive constant $K$ such that for any
two points $x$ and $y$ in $M$ and $z=(x+y)/2$,
\[ K^{-1} \leq \frac{|h(x)-h(z)|}{|h(z)-h(y)|} \leq K.\]
We call the smallest such constant $QC_{h}$ the quasisymmetrical constant of $h$.
We say two mappings $f$ and $g$ from $M$ to itself are
quasisymmetrically
conjugate if they are topologically conjugate and the conjugating
homeomorphism is quasisymmetrical.
In this subsection, we study the quasisymmetrical property of a
conjugating homeomorphism between two geometrically finite
one-dimensional mappings.

\vskip5pt
\noindent {\bf \S 2.4.1 Bounded nearby geometry.}

\vskip3pt
The bounded geometry is a nice geometric property of a
hierarchical structure of intervals. But it is still not enough to get
the quasisymmetrical property of
the conjugating mapping. So we introduce
another concept, bounded nearby geometry.

\vskip3pt
{\sc Definition 4.}
{\em We say a sequence $\eta =\{ \eta_{n}\}_{n=1}^{\infty}$ of nested
partitions of $M$ has bounded nearby geometry if there is a positive
constant $K$ such that for any pair $J$ and $I$ in $\eta_{n}$ with a
common endpoint, the ratio
$|J|/|I| \geq K$. We call the biggest such constant $NC_{f}$ bounded
nearby geometry constant.}

\vskip5pt
{\sc Lemma 2.}
{\em Suppose $f$ is geometrically finite and $\eta =\{ \eta_{n}\}_{n=1}^{\infty}$ is the induced
sequence of nested partitions of $M$ by $f$. Then $\eta$ has bounded
nearby geometry.}

\vskip3pt
{\it Proof.}
We use the same notations as that in the proof of Theorem A. Let
$n_{2}>n_{0}$ be a positive integer such that if two intervals $I$ and $J$
in $\eta_{n_{2}}$ with a
common endpoint, then either both of them are in ${\cal U}$ or both of them
are not in ${\cal U}_{1}$, where
${\cal U}_{1}$ is the union of the critical intervals in $\eta_{n_{2}}$. Let
$K_{1}>0$ be the minimum of ratios, $|J|/|I|$, where $J$ and $I$ are intervals in $\eta_{n_{2}}$ with a
common endpoint.

\vskip5pt
For any $n>n_{2}$ and any two intervals $J$ and $I$ in $\eta_{n}$ with a
common endpoint,
let $J_{i}=f^{\circ i}(J)$ and $I_{i}=f^{\circ i}(I)$ for $i=0$, $\cdots$, $n-n_{2}$.  We consider the intervals $\{ J_{i}\}_{n=0}^{n-n_{2}}$
and $\{ I_{i}\}_{n=0}^{n-n_{0}}$ in the two cases. One is that for some
$0< l\leq n-n_{0}$, $J_{l}=I_{l}$. The other is that $J_{i}$ and $I_{i}$ are different (but they have a common endpoint) for every $i$.

\vskip5pt
For the first case, let $l$ be the smallest such integer, then the common endpoint
of $J_{l+1}$ and $I_{l+1}$ is an extremal critical point of $f$ (this
means that it is either maximal or minimal point of $f$). It
is easy to see now that there is a positive constant $K_{2}$ such that
\[ \frac{|J_{l+1}|}{|I_{l+1}|} \geq K_{2}.\]
Now we use the arguments like that of the second case in the proof of Theorem
A to verify that there is a positive constant $K_{3}$ such that
\[ \frac{|J|}{|I|} \geq K_{4}=K_{3}K_{2}.\]

\vskip3pt
For the second case, again use the arguments like that of the second case in
the proof of Theorem A to demonstrate that there is a positive constant
$K_{5}$ such that
\[ \frac{|J|}{|I|} \geq K_{6}=K_{5}K_{1}.\]
The bounded nearby geometry constant $NC_{f}$ is greater than the maximum of
$K_{4}$ and $K_{6}$.

\vskip5pt
\noindent {\bf \S 2.4.2 Quasisymmetry.}

\vskip3pt
One of the consequences of these properties, bounded geometry and
Bounded nearby geometry is the quasisymmetrical classification
of geometrically finite one-dimensional mappings as follows.

\vskip3pt
{\sc Theorem B.} {\em Suppose $f$ and $g$ are geometrically finite
and topologically conjugate. They are then quasisymmetrically conjugate.}

\vskip3pt
{\it Proof.} Suppose $h$ is the topological conjugacy between $f$ and $g$
and $h\circ f=g\circ h$. Suppose $BC_{f}$, $NC_{f}$, $BC_{g}$
and $NC_{g}$ are the bounded geometry constants and bounded nearby geometry
constants of the induced sequences $\{ \eta_{n,f}\}_{n=1}^{\infty}$ and $\{
\eta_{n,g}\}_{n=1}^{\infty}$ of nested partitions of $M$ by $f$
and $g$ and $\lambda_{n, f}$ and
$\lambda_{n, g}$ are the maximum lengthes of the intervals in $\eta_{n, f}$
and $\eta_{n, g}$, respectively.

\vskip3pt
For any $x < y$ in $M$, let $z$ be the midpoint $(x+y)/2$ of them.
Suppose $N>0$ is the smallest integer such that there is an interval $I$ in
$\eta_{N}$ and is contained in $[x,y]$ (see Figure 1, 2 and 3).

\vskip20pt
\centerline{\psfig{figure=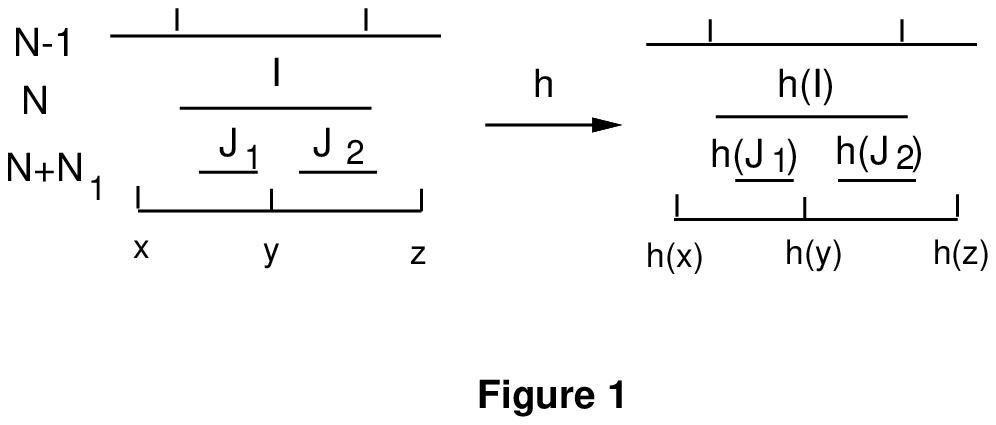}}

\vskip20pt
\centerline{\psfig{figure=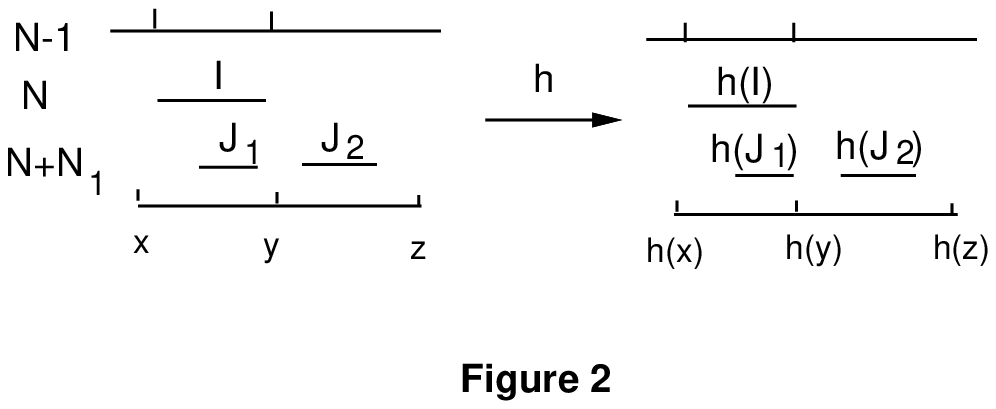}}

\vskip10pt
\centerline{\psfig{figure=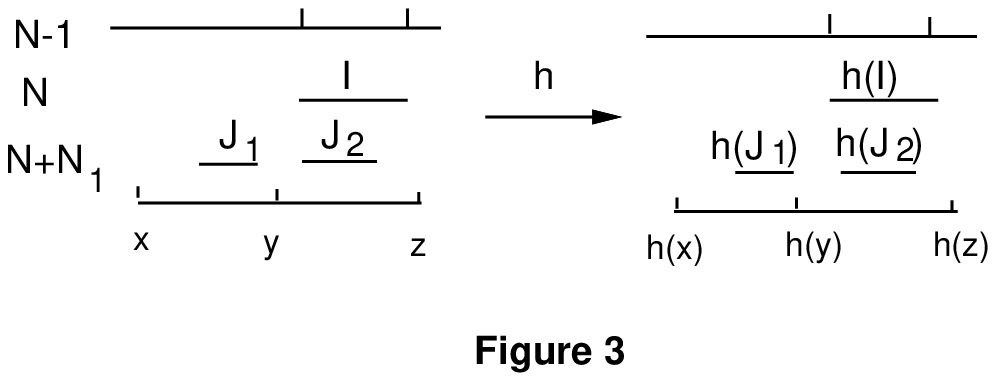}}

\vskip20pt
Because $f$ and $g$ are both geometrically finite, the $n^{th}$-partitions  
$\eta_{n,f}$ and $\eta_{n,g}$ tend to zero exponentially with $n$. We can
find two constants $L_{f}=K(BC_{f}, NC_{f})>0$ and $0< \mu_{f}=\mu (BC_{f},
NC_{f})<1$ such that $\lambda_{n, f} \leq L_{f} (\mu_{f})^{n}$ for any $n>0$.
Moreover, we can
find a positive integer $N_{1}=N_{1}(BC_{f}, NC_{f})$ such
that there are intervals $J_{1}$ and $J_{2}$ in $\eta_{N+N_{1}}$
contained in $[x,z]$
and $[z, y]$, respectively.
By the bounded geometry and bounded nearby geometry, we can find a
constant $K=K(N_{1}, BC_{g}, NC_{g})$ (see Figure 1, 2 and 3)
such that
\[ K^{-1} \leq \frac{|h(x)-h(z)|}{|h(z)-h(y)|} \leq K.\]
The quasisymmetric constant $QC_{h}$ is less than  $K$.

\vskip5pt
\noindent {\bf \S 2.5 Closeness under quasisymmetrical conjugacy.}

\vskip3pt
Another consequence of these properties, bounded
geometry and bounded nearby geometry, is that
geometrically
finite one-dimensional mappings is closed under quasisymmetrical
conjugacy in the space of $C^{1+\alpha}$-mappings with only power law
critical points as follows.

\vskip5pt
{\sc Theorem C.} {\em If a $C^{1+\alpha }$-mapping
$f:M\mapsto M$ for some $0< \alpha \leq 1$ with only power law critical
points is conjugate to a geometrically
finite one-dimensional mapping, then it is also a geometrically finite
one-dimensional mapping.}

\vskip5pt
{\it Proof.} The proof of this theorem is the use of the
quasisymmetrical property of the conjugating homeomorphism.

\vskip10pt
\centerline {\Large \bf \S 3 Examples Of Geometrically Finite}

\centerline{\Large \bf One-dimensional Mappings}

\vskip5pt
The definition of a geometrically finite one-dimensional mapping is
quit abstract. To concrete it, we show some examples. The
main theorem in this section is Theorem D.

\vskip5pt
\noindent {\bf \S 3.1 A $C^{3}$-mapping with nonpositive Schwarzian
derivative.}

\vskip5pt
Suppose $f: M\mapsto M$ is a $C^{3}$-mapping.  The Schwarzian derivative of $f$ is defined by
\[ S(f)(x) = \frac{f'''(x)}{f'(x)} -\frac{3}{2} \Big(
\frac{f''(x)}{f'(x)}\Big)^{2}.\]
We say $f$ has nonpositive Schwarzian derivative if $S(f)(x)\leq 0$ for all $x$ in $M$
and
has nonnegative Schwarzian derivative if
$S(f)(x)\geq 0$ for all $x$ in $M$. We note that a
$C^{3}$-diffeomorphism $f$
has nonpositive Schwarzian derivative if and only if the inverse of $f$
has nonnegative Schwarzian derivative.
The first example of a geometrically finite one-dimensional mapping is the
following:

\vskip5pt
{\sc Example 1.} {\em A $C^{3}$-mapping $f: M\mapsto M$ with finitely
many,
critically finite, nonperiodic power law critical points and nonpositive
Schwarzian derivative.}

\vskip5pt
Suppose $I$ and $J$ are two intervals and $g$ is a $ C^{3}$-diffeomorphism
from $I$ to $J$. A measure of the nonlinearity of $g$ is the
function $n(g) = g^{\prime \prime }/g^{\prime }$. If the absolute
value of $n(g)$ on $I$ is bounded above by a positive constant $C$,
then the distortion $|g'(x)|/|g'(y)|$ of $g$ at any pair $x$ and $y$
in $I$ is bounded above by $exp(C|x-y|)$.
Suppose $d_{I}(x)$
is the distance from $x$ to the
boundary of $I$.

\vskip5pt
{\sc Lemma 3} (the $C^{3}$-Koebe distortion lemma).
{\em Suppose $g$ has nonnegative Schwarzian derivative. Then
$|n(g)(x)|$ is bounded above by $2/d_{I}(x)$ for any $x$ in $I$.}

\vskip5pt
{\it Proof.} See, for example, [J1] for a proof.

\vskip5pt
{\sc Lemma 4.} {\em
Suppose $f$ is the mapping in Example 1 and
$\eta =\{ \eta_{n} \}_{n=1}^{\infty}$ is the induced sequence
of nested partitions of $M$ by $f$. Then $\eta$ has bounded geometry.}

\vskip5pt
{\it Proof.} The proof is similar to that of Theorem A.  Here we use the
Lemma 3 (the $C^{3}$-Koebe distortion lemma) to replace the role of
Lemma 1 (the $C^{1+\alpha }$-Denjoy-Koebe distortion lemma)
in the proof of Theorem A.

\vskip5pt
{\sc Corollary 1.} {\em
Suppose $f$ is the mapping in Example 1 and
$\eta =\{ \eta_{n} \}_{n=1}^{\infty}$ is the induced sequence
of nested partitions of $M$ by $f$.
Then the $n^{th}$-partition $\eta_{n}$ induced by $f$ goes to zero
exponentially with $n$.}

\vskip5pt
{\it Proof.} Suppose $l_{1}$ is the number of the intervals in the first
partition $\eta_{1}$.  Because every critical point of $f$ is not periodic and
critically finite and every periodic point of $f$ is
expanding,
we can find an integer $k>0$ such that every interval
in the first partition $\eta_{1}$ contains at least two but no more than
$kl_{1}$ intervals in the $k^{th}$-partition. By using the bounded geometry,
we can prove this corollary.

\vskip5pt
That Example 1 is geometrically finite follows
from Lemma 3, 4 and Corollary 1.

\vskip10pt
\noindent {\bf \S 3.2 A $C^{1}$-mapping with bounded variation
derivative.}

\vskip5pt
We say a function $u:U\mapsto {\bf R}^{1}$ has bounded variation if
\[ Var(u) = \sup_{x_{1}<\cdots < x_{l}\in U} \sum_{i=1}^{l-1}
 |u(x_{i})-u(x_{i+1})| < +\infty  \]
where $U$ is a subset of $M$.

\vskip5pt
We definite the term $C^{1+bv}$. Suppose $f: M\mapsto M$ is a $C^{1}$-mapping with only
power law critical points. We use $CP=\{ c_{1}$, $\cdots $, $c_{d}\}$
to denote the set of critical points of $f$ and use
$\Gamma =\{ \gamma_{1},$ $\cdots $, $\gamma_{d} \}$ to denote the
the corresponding exponents. Suppose $\eta_{0}$ is the set of intervals
in the complement of the set $CP$ of critical points of $f$.

\vskip5pt
{\sc Definition 5.} {\em We say the mapping $f$ is a $C^{1+bv }$-mapping if

\vskip3pt
$(i)$ the restrictions of $f$ to the intervals in $\eta_{0}$ are
$C^{1}$ with bounded variation derivatives and

\vskip5pt
$(ii)$ for every critical point $c_{i}$ of $f$, there
is a small neighborhood $U_{i}$ of $c_{i}$ such that the functions
$r_{-,i}(x)=f^{\prime}(x)/|x-c_{i}|^{\gamma_{i}-1}$,
$x<c_{i}$ and $x_{i}\in U_{i}$, and
$r_{+,i}(x)=f^{\prime}(x)/|x-c_{i}|^{\gamma_{i}-1}$, $x>c_{i}$ and
$x_{i}\in U_{i}$, have bounded variations.}

\vskip5pt
The main theorem in this section is the following.

\vskip5pt
{\sc Theorem D.} {\em Suppose $f:M\mapsto M$ is a $C^{1+\alpha}$- , for some $0<
\alpha \leq 1$, and
$C^{1+bv }$-mapping with finitely many, critically finite, nonperiodic
power
law critical points and only expanding periodic points and suppose $\eta
=\{
\eta_{n} \}_{n=1}^{\infty}$ is the induced sequence of nested partitions
of $M$ by $f$. Then $\eta$ has bounded geometry.}

\vskip5pt
We prove this theorem by several lemmas. Suppose $f$ is the mapping in Theorem
D. We say an interval $I$ is a $n$-homterval of $f$ if the
restriction of the $i^{th}$-iterate of $f$ to $I$
is a homeomorphism from $I$ to $I_{i}=f^{\circ i}(I)$ for any $i=0$,
$1$, $\cdots$, $n$.  If, moreover, the intervals $\{ I_{i} \}_{i=0}^{n}$
have
pairwise disjoint interiors, then we call it a $n$-wandering homterval.

\vskip5pt
{\sc Lemma D1.} {\em There are
constants $A$, $B>0$ such that for
any $n$-wandering homterval $I$ of $f$ and points $x$ and $y$ in $I$,
\[ \frac{ |(f^{\circ n})^{\prime }(x)|}{|(f^{\circ n})^{\prime}(y)|}
\leq  \exp \Big(A +\frac{B}{D_{x_{n}y_{n}, \partial I_{n}}} \Big) \]
where $x_{n}=f^{\circ n}(x)$, $y_{n}=f^{\circ n}(y)$, $I_{n}=f^{\circ
n}(I)$ and $D_{x_{n}y_{n}, \partial I_{n}}$ is the distance between $\{
x_{n}, y_{n}\}$ and the boundary of $I_{n}$.}

\vskip5pt
{\it Proof.} The idea of the proof of this lemma is the same as that of the proof of the
$C^{1+\alpha }$-Denjoy-Koebe distortion lemma in [J2]. We outline
the proof here.

\vskip5pt
Suppose $U_{i}$ is the set in
Definition 5.
We say an interval in $\eta_{n}$ is an critical interval if one of its
endpoints is a critical point of $f$.  Suppose $n_{0}$ is a
positive integer such that
every critical interval $I$ in $\eta_{n_{0}}$ is contained in
some $U_{i}$ and one of its endpoints
is not in the critical orbits $\cup_{n=0}^{\infty}f^{\circ
n}(CP)$.

\vskip5pt
The ratio, $f^{\circ n}(x)/ f^{\circ n}(y)$, equals the product
$\prod_{i=0}^{n-1}f^{\prime}(x_{i})/f^{\prime}(y_{i})$ where
$x_{i}=f^{\circ
i}(x)$ and $y_{i}=f^{\circ i}(y)$.  We divide this product into two products,
\[ \prod_{x_{i}, y_{i}\in {\cal U}}\frac{f'(x_{i})}{f'(y_{i})} \hskip5pt
and \hskip5pt \prod_{x_{i}, y_{i}\in {\cal V}
}\frac{f'(x_{i})}{f'(y_{i})}\]
where ${\cal U}$ stands for the union of all the critical intervals in
$\eta_{n_{0}}$
 and ${\cal V}$ stands for the union of all the noncritical intervals in
$\eta_{n_{0}}$.
The second product is bounded by $\exp (Var(f^{\prime})/\beta )$, where
$\beta >0 $ is
the minimum of the absolute value of the restriction of the derivative $f'$ to ${\cal V}$.

\vskip5pt
Let
$\tilde{r}_{-,i}(x)=|f(x)-f(c_{i})|/|x-c_{i}|^{\gamma_{i}}$,
$x<c_{i}$ and $x_{i}\in U_{i}$, and
$\tilde{r}_{+,i}(x)=|f(x)-f(c_{i})|/|x-c_{i}|^{\gamma_{i}}$,
$x>c_{i}$ and $x_{i}\in U_{i}$. Then both of them have bounded
variations. We may write the first product into
\[  \prod_{x_{i}, y_{i}\in {\cal U}}\frac{|f'(x_{i})|}{|f'(y_{i})|} =
\prod_{x_{i}, y_{i}\in
{\cal U}}\frac{|r_{a_{i}b_{i}}(x_{i})|}{|r_{a_{i}b_{i}}(y_{i})|}
\frac{(\tilde{r}_{a_{i}b_{i}}(y_{i}))^{\gamma_{i}-1}}{(\tilde{r}_{a_{i}b_{i}}(x_{i}))^{\gamma_{i}-1}}
\frac{|f(x_{i})-f(c_{b_{i}})|^{m_{b_{i}}}}{|f(y_{i})-f(c_{b_{i}})|^{m_{b_{i}}}}\]
where $a_{i}$ is $+$ or $-$, $b_{i}$
the integer such that $x_{i}$ and $y_{i}$ are in $U_{b_{i}}$ and
$m_{b_{i}}=1-1/\gamma_{b_{i}}$.
The first two products satisfy that
\[ \prod_{x_{i}, y_{i}\in
{\cal U}}\frac{|r_{a_{i}b_{i}}(x_{i})|}{|r_{a_{i}b_{i}}(y_{i})|}
\frac{(\tilde{r}_{a_{i}b_{i}}(y_{i}))^{\gamma_{i}-1}}
{(\tilde{r}_{a_{i}b_{i}}(x_{i}))^{\gamma_{i}-1}}
\leq \]
\[ \exp \Big( \sum_{i=1}^{d_{1}} \Big( Var(r_{i+}) +Var(r_{i-}) +
\frac{1}{\gamma_{i}-1}\Big( Var(\tilde{r}_{i+})
+Var(\tilde{r}_{i-}) \Big) \Big) \Big) .\]

\vskip5pt
To estimate the last product, we write each
\[ \frac{f(x_{i})-f(c_{b_{i}})}{f(y_{i})-f(c_{b_{i}})}=1+
\frac{f(x_{i})-f(y_{i})}{f(y_{i})-f(c_{b_{i}})}\]
for $x_{i}$ and $y_{i}$ in ${\cal U}$
and
\[ \prod_{x_{i}, y_{i}\in {\cal U}} \frac{|f(x_{i})-f(c_{b_{i}})|^{m_{b_{i}}}}{|f(x_{i})-f(c_{b_{i}})|^{m_{b_{i}}}}
\leq  \exp \Big( \frac{1}{m_{i}}\sum_{k=1}^{l} \log \Big( 1+
\frac{|f(x_{i_{k}})-f(y_{i})|}{|f(x_{i_{k}})-f(c_{b_{i_{k}}})|} \Big) \Big) \]
where $i_{1} < \cdots < i_{l}$.

\vskip5pt
Because each critical point of $f$ is mapped eventually to an expanding
periodic
point and all the periodic points of $f$ are expanding,
there is a positive constant $K_{1}$ (by the naive distortion lemma in [J2]) such
that
\[ \frac{|f(x_{i_{1}})-f(y_{i_{1}})|}{|f(x_{i_{1}})-f(c_{b_{i_{1}}})|}
 \leq K_{1}\frac{|x_{n}-y_{n}|}{D_{x_{n}y_{n}, \partial I_{n}}},\]
and
\[ \frac{|f(x_{i_{k}})-f(y_{i_{k}})|}{|f(x_{i_{k}})-f(c_{b_{i_{k}}})|} \leq
K_{1}\frac{|x_{i_{k-1}}-y_{i_{k-1}}|}{|y_{i_{k-1}}-f^{i_{k}-i_{k-1}}(c_{b_{i_{k}}})|} \]
for any $0< k \leq l$. We may assume that
$f^{i_{k}-i_{k-1}-1}(c_{b_{i_{k}}})$ is not a critical point of
$f$.
Otherwise, $f^{i_{k}-i_{k-1}}(c_{i_{k}})=c_{i_{k-1}}$ and note that
there are only finitely many critical chains of $f$ like that in the
proof of Theorem A. Let $L$ be
the minimum of lengths of the critical intervals in $\eta_{n_{0}}$. Then
\[ \frac{|f(x_{i_{k}})-f(y_{i_{k}})|}{|f(x_{i_{k}})-f(c_{b_{i_{k}}})|}
\leq K_{1}\frac{|x_{i_{k-1}}-y_{i_{k-1}}|}{L}, \]
and moreover, there are constants $K_{2}$, $K_{3} >0$ such that
\[ \prod_{x_{i}, y_{i}\in {\cal U}} \frac{|f(x_{i})-f(c_{b_{i}})|^{m_{b_{i}}}}{|f(x_{i})-f(c_{b_{i}})|^{m_{b_{i}}}}
\leq  K_{2}+\frac{K_{3}}{D_{x_{n}y_{n}, \partial I_{n}}}.\]

\vskip5pt
Combining all the estimates together, we get two positive constants $A$ and $B$.

\vskip5pt
{\sc Lemma D2.} {\em Every $\infty$-homterval $I$ of $f$ is
an $\infty$-wandering homterval.}

\vskip5pt
{\it Proof.} Suppose there are integers $m>n>0$ such that $I_{n}$ and
$I_{m}$
are overlap. Let $k=n-m$, then $I_{0}$ and $I_{k}$ are overlap,
and moreover, $I_{lk}$ and $I_{(l+1)k}$ are overlap for any $l>0$. Let
$T=\cup_{l=0}^{\infty}I_{kl}$.	It is a connected interval of $M$ and
$f^{\circ k}: T\mapsto T$ is a homeomorphism.  Then $f^{\circ k}$ has to
have a fixed point which is not topologically expanding. This
contradiction proves the lemma.

\vskip5pt
From Lemma D1 and Lemma D2, we have the following lemma.

\vskip5pt
{\sc Lemma D3.} {\em The maximal length of the intervals in $\eta_{n}$
tends to zero as $n$ goes to infinity.}

\vskip5pt
{\it Proof.}
Suppose there is an $\epsilon_{0}>0$ such that for any positive integer $n$,
there is an interval $I_{n}\in \eta_{n}$ with $|I_{n}|>\epsilon_{0}$.
Because $M$ is a compact manifold, there is a subset $\{
n_{i}\}_{i=1}^{\infty}$
of the integers such that $I_{n_{i}}$ goes to an interval
$\tilde{I}$ as $i$ goes to infinity and the length of $\tilde{I}$
is greater than $\epsilon_{0}$.  There is an interval
    $I\subset \tilde{I}$ such that $I\subset I_{n_{i}}$ for large $i$.
    The restriction of the $i^{th}$-iterate of $f$ to $I_{n_{i}}$ is
   an embedding for any $i\leq n_{i}$. Hence $I$ is an
$\infty$-homterval of
    $f$, and moreover, it is an $\infty$-wandering homterval. Suppose $I$ is
    a maximal such interval. Let $T_{n}\supset I$ be the maximal
    $n$-homterval. Then it is again a $n$-wandering homterval. Let
    $L_{n}$ and $R_{n}$ be the intervals in the complement
    of $I$ in $T_{n}$. The lengths of $L_{n}$ and $R_{n}$ go to zero as $n$
tends to infinity. The boundary of $f^{\circ n}(T_{n})$ is contained in
the union of the boundary of $M$ and the set of critical values $f(CP)$
of $f$ for $T_{n}$ is a maximal
    $n$-homterval of $f$. Suppose $\{n_{i}\}_{i=0}^{\infty}$ is a
    subsequence of the integers such that
    the boundary of $f^{\circ n_{i}}(T_{n_{i}})$ are the same for all
    $i$.  By using
    Lemma D1, one of the lengths of $f^{\circ n_{i}}(I\cup L_{n_{i}})$
and
    $f^{\circ n_{i}}(R_{n_{i}}\cup I)$, say
    $f^{\circ n_{i}}(I\cup L_{n_{i}})$,
    has to go to zero as $i$ tends
    to infinity.
    Because every critical point is mapped to a periodic point eventually,
    the interval
    $f^{\circ n_{i}}(I\cup L_{n_{i}})$ tends to a periodic orbit
    eventually. This periodic point is not topologically expanding. The
    contradiction proves the lemma.

\vskip5pt
Recall that in the proof of Lemma D1, ${\cal U}$ stands for the  union of
all the critical intervals in $\eta_{n_{0}}$ and ${\cal V}$ stands for the
union of all the noncritical intervals in $\eta_{n_{0}}$,
where $n_{0}$ is a fixed
positive integer such that every critical interval is contained in $U_{i}$
in Definition 5 and one of its endpoints is not in the critical
orbits
$\cup_{n=0}^{\infty}f^{\circ n}(CP)$.

\vskip5pt
Lemma D4 and Lemma D5 are two of the key lemmas in the proof of Theorem D.

\vskip5pt
{\sc Lemma D4.} {\em There is a constant $K>0$ such that
for an interval $I\in \eta_{n+n_{0}}$, if $I_{i}=f^{\circ i}(I)$ is
in ${\cal V}$ for every $1\leq i\leq n$, then
\[ \frac{|(f^{\circ n})'(x)|}{|(f^{\circ n})'(y)|} \leq K\]
for any $x$ and $y$ in $I$.}

\vskip5pt
{\it Proof.}
If $\{ I_{i}\}_{i=0}^{n-1}$ have pairwise disjoint
interiors, then
\[ |(f^{\circ n})'(x)| \geq \exp
\Big( \frac{Var(f')}{\beta }\Big) \frac{|I_{n}|}{|I|}\]
for any $x\in I$ where $I_{n}=f^{\circ n}(I)\in \eta_{n_{0}}$ and $\beta >0$ is the
minimum of the absolute value of $f'|{\cal V}$.
By using this fact and Lemma D3,
we can find a constant $\nu >1$ such that for a periodic point $p$ of
$f$, if $p_{i}=f^{\circ i}(p)$ is in ${\cal V}$ for every  $i\geq 0$, then the eigenvalue $|(f^{\circ k})'(p)|\geq
\nu$ where $k$ is the period of $p$.

\vskip5pt
By the naive distortion lemma (see [J1] or [J2]), we have that
\[ \frac{|(f^{\circ n})'(x)|}{|(f^{\circ n})'(y)|} \leq
\exp \Big( \frac{K_{1}}{\beta }\sum_{i=0}^{n-1}|I_{i}|^{\alpha } \Big) \]
for any $x$ and $y$ in $I$ where $K_{1}$ is a positive constant and $c$ is
the minimum of the absolute value of $f'|{\cal V}$.

\vskip5pt
Suppose $I_{0}$, $\cdots $, $I_{k-1}$ have pairwise disjoint interiors and
$I_{k}\subset I_{0}$.  There is a periodic point $p$ of period $k$ in
$I_{0}$. Again by using the naive distortion lemma, there is a constant
$K_{2}>0$ such that
\[ |I_{lk+i}|\leq \frac{K_{2}}{\nu^{l}} |I_{i}|\]
for all $l > 0$ and $0\leq i <k$ where $K_{2}>0$ is a constant. Last two
inequalities imply Lemma D4.

\vskip5pt
We say a critical point of $f$ is pure if
it is not in the post-critical orbits
$\cup_{n=1}^{\infty}f^{\circ n}(CP)$. We say an interval $I$ is a
pure critical interval in $\eta_{n_{0}}$ if one of its endpoint
is pure
critical point. Remember that the other endpoint of $I$ is not in the
critical orbits $\cup_{n=0}^{\infty} f^{\circ n}(CP)$.

\vskip5pt
{\sc Lemma D5.} {\em There is a constant $K>0$ such that
for an interval $I\in \eta_{n+n_{0}}$, if $I_{n}=f^{\circ n}(I)$ is in
a pure critical interval in $\eta_{n_{0}}$, then
\[ \frac{|(f^{\circ n})'(x)|}{|(f^{\circ n})'(y)|} \leq K\]
for any $x$ and $y$ in $I$.}

\vskip5pt
{\it Proof.} By the similar arguments to the proof of Lemma D1
and that $I_{n}$ is far to the post-critical orbit
$\cup_{n=1}^{\infty}f^{\circ n}(CP)$, we can find a
positive constant $K_{1}$ such that
if $\{ I_{i}\}_{i=0}^{\infty}$ have pairwise disjoint interiors, then
\[ |(f^{\circ n})'(x)| \geq \exp \Big(
K_{1} \Big) \frac{|I_{n}|}{|I|}\] for any $x\in I$.  Using
this fact and Lemma
D3,
we can find a constant $\nu >1$ such that for any periodic point $p$
in a pure critical interval in
$\eta_{n_{0}}$, the eigenvalue $|(f^{\circ k})'(p)|\geq \nu$ where $k$
is the period of $p$.

\vskip5pt
By the version of the $C^{1+\alpha}$-Denjoy-Koebe distortion lemma in [J2], there is
a constant $K_{2}>0$ such that
\[ \frac{|(f^{\circ n})'(x)|}{|(f^{\circ n})'(x)|} \leq
\exp \Big( K_{2} \sum_{i=0}^{n-1}|I_{i}|^{\alpha} \Big) \]
for any $x$ and $y$ in $I$ where $K_{2}$ is a positive constant.

\vskip5pt
Suppose $I_{0}$, $\cdots $, $I_{k-1}$ have pairwise disjoint interiors and
$I_{k}\subset I_{0}$.  There is a periodic point $p$ of period $k$ in
${I_{0}}$. Again by Lemma 1
and the naive distortion lemma (see [J1] or [J2]),
there is a constant $K_{3}>0$
such that
\[ |I_{lk+i}|\leq \frac{K_{3}}{\nu^{l}} |I_{i}|\]
for all $l > 0$ and $0\leq i <k$. The last two inequalities imply Lemma D5.

\vskip5pt
{\it Proof of Theorem D.} The proof of Theorem D is now similar to the
proof of Theorem A. Here we use Lemma D4 to replace the role of the
naive distortion lemma and use Lemma D5 to replace the role of Lemma 1 (the $
C^{1+\alpha }$-Denjoy-Koebe distortion lemma).

\vskip5pt
{\sc Corollary D1.} {\em The maximum $\lambda_{n}$ of lengths of the
intervals
in $\eta_{n}$ tends to zero exponentially with $n$.}

\vskip5pt
Theorem D and Corollary D1 provide another example of a
geometrically finite one-dimensional mapping.

\vskip5pt
{\sc Example 2.} {\em A $C^{1+\alpha}$-, for some $0< \alpha \leq 1$, and
$C^{1+bv}$-mapping $f:M\mapsto M$
 with finitely many, critically finite, nonperiodic power law critical
points and only expanding periodic points.}

\vskip5pt
In Theorem D and in Example 2, the condition that $f$ is a
$C^{1+\alpha}$-, for some $0<\alpha \leq 1$, and $C^{1+bv }$-mapping can not be
weakened
to the condition that $f$ is a $C^{1+\alpha }$-mapping for there is a
counterexample in [J1].
The construction
of the counterexample in [J1] is like the
construction of the Denjoy counterexample in circle diffeomorphisms and
this example is not topologically conjugate to any geometrically
finite one-dimensional mapping.

\vskip5pt
\noindent {\bf \S 3.3 A question on $C^{1+\alpha}$-mappings with
expanding periodic points.}

\vskip3pt
In Theorem C,
the conditions that a $C^{1+\alpha }$-mapping for some $0< \alpha \leq 1$
with only power law critical points is
quasisymmetrically
conjugate to a geometrically finite one-dimensional mapping can not be
weakened
to the condition that a $C^{1+\alpha }$-mapping for some $0<\alpha \leq 1$
with only power law critical point is topologically
conjugate to
a geometrically finite one-mapping for there is an easy counterexample
$f: [-1, 1]\mapsto [-1,1]$ with the neutral fixed point $-1$, namely
$f'(-1)=1$ (see Figure 4).

\vskip10pt
\centerline{\psfig{figure=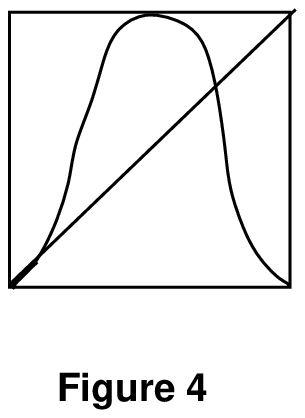}}

\vskip10pt
The graph in Figure 4 suggests a question
as follows.

\vskip5pt
{\sc Question 1.} {\em Suppose $f: M\mapsto M$ is a
$C^{1+\alpha}$-mapping for some $0< \alpha \leq 1$
with
only power law critical points
and only expanding periodic points and
is
topologically
conjugate to a geometrically finite one-dimensional mapping.
Is $f$
geometrically finite ?}

\vskip5pt
The answer of this question may be negative.
But we do not have a concrete counterexample yet.
The reader may refer to the construction of the counterexample in [J1]
and Lemma D4 and Lemma D5 in this paper.


\vskip20pt
\begin{thebibliography}{99}

\bibitem[A]{a} {\sl L. A. Ahlfors,} [1966]. Lectures on Quasiconformal
Mappings. D. von Nostrand Company, Inc.

\bibitem[BL]{bl} {\sl A. M. Blokh and M. Yu. Lyubich}, [1990].
Measurable Dynamics of S-Unimodal Maps of the Interval. {Preprint, IMS,
SUNY at Stony Brook.}

\bibitem[H]{he} {\sl M. R. Herman}, [1989]. Seminars given by M. Herman
in IAS.

\bibitem[Ja]{ja} {\sl M. Jakobson}, [1989]. Quasisymmetric Conjugacy
for
Some One-dimensional Maps Inducing Expansion. {\em preprint.}

\bibitem[J1]{yj1} {\sl Y. Jiang}, [1990]. Generalized Ulam-von Neumann
Transformation. {\em Thesis, Graduate School of CUNY}.

\bibitem[J2]{ji2} {\sl Y. Jiang}, [1990]. Dynamics of Certain Smooth
One-dimensional Mappings, I. The $C^{1+\alpha }$-Denjoy-Koebe distortion lemma.
{\em Preprint}.

\bibitem[J3]{ji3} {\sl Y. Jiang}. Dynamics of Certain Smooth
One-dimensional Mappings, III. Scaling function geometry.
{\em Preprint}.

\bibitem[M]{ma} {\sl R. Ma\~{n}e}. Hyperbolicity,
Sinks
and Measure in One Dimensional Dynamics, [1985]. {\em Comm. in Math.
Phys. {\bf 100}, 495-524} and {\em Erratum Comm. in Math. Phys. {\bf
112}, (1987) 721-724}.

\bibitem[Mi]{mis} {\sl M. Misiurewicz}. Absolutely
Continuous Measures for Certain Maps of An Interval, [1978].
{\em Inst. Hautes. \'Etudes. Sci. Publ. Math. 1978, No. {\bf 53}, 17 -
51.}

\bibitem[MS]{ms} {\sl W. de Melo and S. van Strien}. A Structure Theorem
in One Dimensional Dynamics, [1986].
{\em Report 86-29, Delft University of Technology, Delft}.

\bibitem[MT]{mt} J. Milnor and W. Thurston, [1977]. On iterated
maps of the interval I and II. {\em Preprint, Princeton University
Press:Princeton.}

\bibitem[S1]{ds} {\sl D. Sullivan}. Class Notes, [1989].

\bibitem[S2]{ds} {\sl D. Sullivan}. Bounded Structure of
Infinitely Renormalizable Mappings, [1989]. {\em Universality in Chaos,
$2^{n}$ Edition, Adam Hilger, Bristol, England.}.

\bibitem[SW]{sw} {\sl G. Swiatek}, [1990].
Private Conversations.

\bibitem[Y]{y} {\sl J. Yoccoz}, [1989]. Seminars given by M. Herman in
IAS.
\end{thebibliography}
\end{document}